\documentclass[graybox]{svmult}

\usepackage{type1cm}        %
\usepackage{makeidx}         %
\usepackage{graphicx}        %
\usepackage{multicol}        %
\usepackage[bottom]{footmisc}%
\usepackage{listings}
\usepackage{tikz}
\usepackage{pgfplots}
\usepackage{siunitx}
\pgfplotsset{compat=1.18}
\usepackage{caption,subcaption}
\usepackage[dvipsnames]{xcolor}
\usepackage{nicefrac}
\usepackage{microtype}
\usepackage{algorithm2e}

\usepackage[colorlinks]{hyperref}

\usepackage{newtxtext}       %
\usepackage[varvw]{newtxmath}       %
\usepackage{todonotes}
\usepackage[normalem]{ulem}

\makeindex             %

\begin{document}

\title*{Low-Rank Solvers for Energy-Conserving Hamiltonian Boundary Value Methods}
\titlerunning{Low-Rank Solvers for HBVMs}
\author{Fabio Durastante\orcidID{0000-0002-1412-8289} and\\ Mariarosa Mazza\orcidID{0000-0002-8505-6788}}
\authorrunning{F. Durastante, M. Mazza}
\institute{Fabio Durastante \at University of Pisa, Largo Bruno Pontecorvo, 5, Pisa, 56127, \email{fabio.durastante@unipi.it}
\and Mariarosa Mazza \at University of Rome ``Tor Vergata'', via della Ricerca Scientifica, 1, Roma, 00133, \email{mariarosa.mazza@uniroma2.it}}
\maketitle

\abstract*{We study energy-conserving Hamiltonian Boundary Value Methods (HBVMs) for Hamiltonian systems, which arise in applications where long-term preservation of energy and symplecticity is essential. HBVMs are multi-stage schemes whose stage equations reformulate as matrix equations with a low-rank right-hand side. For linear systems, we exploit this structure directly via Krylov projection solvers. For nonlinear systems, we leverage it within simplified Newton iterations and as a preconditioner in a Newton--Krylov framework, combined with adaptive time-stepping for robust convergence. Numerical experiments on semi-discretized wave equations demonstrate the efficiency and robustness of the proposed approach.}

\abstract{We study energy-conserving Hamiltonian Boundary Value Methods (HBVMs) for Hamiltonian systems, which arise in applications where long-term preservation of energy and symplecticity is essential. HBVMs are multi-stage schemes whose stage equations reformulate as matrix equations with a low-rank right-hand side. For linear systems, we exploit this structure directly via Krylov projection solvers. For nonlinear systems, we leverage it within simplified Newton iterations and as a preconditioner in a Newton--Krylov framework, combined with adaptive time-stepping for robust convergence. Numerical experiments on semi-discretized wave equations demonstrate the efficiency and robustness of the proposed approach.}

\section{Introduction}
\label{sec:1}

In this paper we are interested in the solution of ordinary differential equations (ODEs) of the form
\begin{equation}\label{eq:the_conservative_problem}
\begin{cases}
    \mathbf{y}'(t) = \mathcal{J} \nabla \mathrm{H}(\mathbf{y}(t)), \quad t \in [0,T],\\
    \mathbf{y}(0) = \mathbf{y}_0,
\end{cases} \quad \mathbf{y}_0 \in \mathbb{R}^{2m},
\end{equation}
where 
\[
\mathcal{J} = \begin{bmatrix}
     & I_m \\
   - I_m & 
\end{bmatrix}, \quad \mathrm{H} : \mathbb{R}^{2m} \rightarrow \mathbb{R},
\]
with $I_m$ the $m \times m$ identity matrix. This class of problems, often referred to as Hamiltonian systems, arises naturally in a wide range of applications where conservation laws and symplectic structure are commonly present. Typical examples include the simulation of systems in classical mechanics, such as planetary motion~\cite{10.1093/mnras/stx547} and rigid body dynamics~\cite{Leimkuhler_Reich_2005}, where the Hamiltonian $\mathrm{H}$ represents the total energy. These type of equations also appear in plasma physics~\cite{10.1063/1.1882353}, molecular dynamics~\cite{Leimkuhler1996}, and quantum mechanics~\cite{770c5663-0765-3dbf-a7d6-1543084c96fb}, among others. 

Due to the importance of accurately preserving geometric properties over long time intervals, we focus on the numerical solution of~\eqref{eq:the_conservative_problem} using a specific family of one-step multi-stage methods called Hamiltonian Boundary Value Methods (HBVM) from~\cite{BrugnanoBook,MR2833606,MR2832704} that closely resemble Runge--Kutta (RK) methods. We, therefore, begin by establishing the notation for RK methods applied to ODEs of the form
\begin{equation}\label{eq:theproblem}
\begin{cases}
    \mathbf{y}'(t)=f(\mathbf{y}(t),t), & t \in [0,T],  \\
    \mathbf{y}(0) = \mathbf{y}_0,
\end{cases} \quad \begin{array}{l}
\mathbf{y} \,:[0,T] \rightarrow \mathbb{R}^N,\\
f \,:\,\mathbb{R}^N \times [0,T] \rightarrow \mathbb{R}^N, \quad T > 0 \\
\mathbf{y}_0 \in \mathbb{R}^N,\quad N \in \mathbb{N}.
\end{array}
\end{equation}

An Implicit RK (IRK) method with $s$ stages for integrating~\eqref{eq:theproblem} requires solving a nonlinear system of algebraic equations at each time-step. The update from time $t_n$ to $t_{n+1} = t_n + h$ is given~by
\begin{equation}\label{eq:steps_Runge_Kutta}
\mathbf{y}^{(n+1)} = \mathbf{y}^{(n)} + h \sum_{i=1}^{s} b_i \mathbf{k}_{i}^{(n)}, \quad n = 0, \ldots, n_t-1, \quad h = \frac{T}{n_t}, \quad n_t \in \mathbb{N},
\end{equation}
where the stage values $\mathbf{k}_{i}^{(n)}$ are computed as
\begin{equation}\label{eq:stages_Runge_Kutta}
\mathbf{k}_{i}^{(n)} = f\left(\mathbf{y}^{(n)} + h \sum_{j=1}^{s} a_{i,j} \mathbf{k}_{j}^{(n)},\, t_n + c_i h\right), \qquad i=1,\ldots,s,
\end{equation}
or equivalently 
\begin{equation}\label{eq:stages_Runge_Kutta_v2}
\mathbf{k}_{i}^{(n)} = \mathbf{y}^{(n)} + h \sum_{j=1}^{s} a_{i,j} f\left(\mathbf{k}_{j}^{(n)},\, t_n + c_i h\right), \qquad i=1,\ldots,s.
\end{equation}

The method is fully specified by its coefficients $a_{i,j}$, weights $b_i$, and nodes $c_i$, for $i,j = 1, \ldots, s$. These are conveniently organized in the Butcher tableau:
\begin{equation}\label{eq:butcher_tableau}
\def\arraystretch{1.2}
\begin{array}[b]{c|c}
\mathbf{c} & A\\
\hline
 & \mathbf{b}^\top
\end{array} \;\raisebox{0.8em}{=}\; 
\begin{array}[b]{c|ccc}
c_1 & a_{1,1} & \ldots & a_{1,s}\\
\vdots & \vdots & \ddots & \vdots\\
c_s & a_{s,1} & \ldots & a_{s,s}\\
\hline
 & b_1 & \ldots & b_{s}
\end{array}\raisebox{0.8em}{.}
\end{equation}

The contribution of this paper is to introduce a low-rank technique for the IRK-like formulation of the HBVMs ~\cite{BrugnanoBook,MR2833606,MR2832704} for~\eqref{eq:the_conservative_problem}.  Such approach is in connection with the one introduced in~\cite{durastante2025stageparallelimplicitrungekuttamethods},  and arises from a similar matrix equation perspective. In the Hamiltonian setting, the resulting matrix equation always contains a known low-rank right-hand side. This structure simplifies the procedure of~\cite{durastante2025stageparallelimplicitrungekuttamethods} considerably: for linear Hamiltonians we directly solve a single matrix equation with the low-rank right-hand side per time-step; for nonlinear Hamiltonians we present two alternative choices. We either reduce the procedure to solving a sequence of linear matrix equations by a simplified Newton method~\cite{MR1790038}, or directly tackle the full Newton step via a matrix equation preconditioned Krylov method~\cite{MR1344684}.

This kind of approach is part of a line of research aimed at achieving efficient solution of the---generally nonlinear---equations for the stages in RK-type methods; see, for example,~\cite{Bellen,Bellen2,MR4744993,Houwen1,Houwen3,Houwen2,Iserles,MR1790038,SantoloSVD,MR4735250,SouthworthII,SouthworthI}. Notably, the works~\cite{Houwen1,Houwen3,Houwen2} introduce a diagonal iteration technique for solving the stage equations, which reduces the problem to a set of block-diagonal systems. In a similar vein,~\cite{Iserles} explores the sparsity pattern of the stage matrix, drawing on tools from sparse linear algebra to minimize inter-stage coupling. Finally, the use of preconditioned Krylov methods inside different linearization strategies for the stage equations is addressed in~\cite{MR4744993,MR1790038,SantoloSVD,MR4735250,SouthworthII,SouthworthI}---in more detail, \cite{MR4744993,MR4735250} employ matrix-structure preserving preconditioners, while the others exploit factorizations of the Butcher tableau matrix $A$ like its singular value decomposition \cite{SantoloSVD} and its eigenvalue decomposition \cite{MR1790038,SouthworthII,SouthworthI}.

The rest of the paper is organized as follows. In Section~\ref{sec:conservative_problems}, we introduce the class of numerical methods to which the low-rank strategy can be effectively applied, and we present the general framework of the proposed procedure in detail. Section~\ref{sec:thelinearcase} focuses on the case where the Hamiltonian function \( \mathrm{H} \) in~\eqref{eq:the_conservative_problem} is linear, allowing for a simplified and illustrative implementation of the method. In Section~\ref{sec:thenonlinearcase}, we extend the approach to handle the more general and practically relevant case of a nonlinear Hamiltonian \( \mathrm{H} \), highlighting the necessary modifications and associated computational considerations. In Section~\ref{sec:numerical_examples} we test the algorithm on some benchmark problems. Finally, in Section~\ref{sec:conclusions} we draw conclusions and highlight possible future directions.

\section{A class of structure preserving methods for conservative problems}\label{sec:conservative_problems}

We begin our investigation by considering the numerical solution of problem~\eqref{eq:the_conservative_problem} using a class of energy-preserving methods known as $\operatorname{HBVM}(k,s)$, where $k, s \in \mathbb{N}$ and, typically, $k \geq s$~\cite{BrugnanoBook,MR2833606,MR2832704}. These methods are designed to preserve the Hamiltonian structure of the system and are formulated as single-step, multi-stage schemes that closely resemble RK methods, following the framework introduced in~\eqref{eq:steps_Runge_Kutta}--\eqref{eq:stages_Runge_Kutta}. In this section, we outline the construction of the $\operatorname{HBVM}(k,s)$ method, highlighting its structural similarities to classical RK methods and its ability to exactly preserve polynomial Hamiltonians of a given degree; see~\cite[Theorem 3.4]{BrugnanoBook}.

Let $\{P_\ell\}_{\ell \geq 0}$ be the sequence of scaled and shifted degree $\ell$ Legendre polynomials on the interval $[0,1]$~\cite[P.~27]{Gautschi}, i.e., the sequence of orthonormal polynomials with respect to the $\omega(x) = 1$ measure on the $[0,1]$ interval
\[
\int_{0}^{1} P_p(x) P_q(x)\,\omega(x)\,\mathrm{d}x = \delta_{p,q} = \begin{cases}
    1, & p = q,\\
    0, & p \neq q.
\end{cases}
\]
To build a single-step approximation to~\eqref{eq:the_conservative_problem}, we start by discretizing the interval $[0,T]$ as in~\eqref{eq:steps_Runge_Kutta},
and restrict ourselves to consider the construction of the stages in the first time-step. We introduce the variable $h \in [0,1]$ and expand~\eqref{eq:the_conservative_problem}~as
\[
\mathbf{y}'(\tau h) = \sum_{\ell \geq 0} P_\ell(\tau) \boldsymbol{\phi}_\ell(\mathbf{y}), \quad \boldsymbol{\phi}_\ell\left(\mathbf{y}(ch)\right) = \int_{0}^{1} P_\ell(c) \mathcal{J} \nabla\mathrm{H}(\mathbf{y}(ch))\,\mathrm{d}c.
\]
To arrive at a discrete version we truncate the series after $s$ terms, and fix the number of quadrature nodes to approximate the integrals describing the $\boldsymbol{\phi}_\ell$ to $k$. 
To lighten the notation, let us set
\begin{equation}\label{eq:generic_discretization}
    f(\mathbf{y}) = \mathcal{J} \nabla\mathrm{H}(\mathbf{y}), \quad \sigma' (\tau h) = \sum_{\ell=0}^{s-1} P_\ell(\tau) \widehat{\boldsymbol{\phi}}_\ell(\sigma), 
\end{equation}
where $\sigma'$ is the approximation of $\mathbf{y}'$ obtained by selecting $s$ and the quadrature rule $\widehat{\boldsymbol{\phi}}_\ell$ for the $\boldsymbol{\phi}_\ell$---hence $\sigma$ approximates $\mathbf{y}$ and is such that $\sigma(0) = \mathbf{y}(0) = \mathbf{y}_0$. To obtain the full scheme, we now integrate~\eqref{eq:generic_discretization} on the first step, i.e., for $\tau \in [0,1]$ obtaining
\[
\sigma(\tau h) = \sigma(0) + h \sum_{\ell=0}^{s-1} \int_{0}^{\tau} P_\ell(x)\,\mathrm{d}x\,\widehat{\boldsymbol{\phi}}_\ell(\sigma) = y_0 + h \sum_{\ell=0}^{s-1} \int_{0}^{\tau} P_\ell(x)\,\mathrm{d}x\,\widehat{\boldsymbol{\phi}}_\ell(\sigma).
\]
We now need an explicit expression for the $\{\widehat{\boldsymbol{\phi}}_\ell(\sigma)\}_{\ell=0}^{s-1}$.  
To proceed, we choose the Gauss--Legendre formula with nodes $\{ c_i \}_{i=1}^{k}$ and weights $\{b_i\}_{i=1}^{k}$, hence
\begin{equation}\label{eq:relation-for-the-phi}
    \widehat{\boldsymbol{\phi}}_j = \sum_{i=1}^{k} b_i P_j(c_i) f(\sigma(c_i h)),
\end{equation}
from which we obtain, for $j=0,\ldots,s-1,$
\begin{equation}\label{eq:what-we-have-to-compute}
\begin{split}
\sigma(c_j h) = & \; y_0 + h \sum_{\ell = 0}^{s-1} \int_{0}^{c_j} P_\ell(x)\,\mathrm{d}x \sum_{i = 1}^{k} b_i P_\ell(c_i) f(\sigma(c_i h)) \\
= & \; y_0 + h \sum_{i = 1}^k b_i \left[ \sum_{\ell = 0}^{s-1} \int_{0}^{c_j} P_\ell(x)\,\mathrm{d}x \cdot P_\ell(c_i) \right] f(\sigma(c_i h)).
\end{split}
\end{equation}
We define the $\mathcal{W}_s \in \mathbb{R}^{k \times s}$ transformation matrix
    \begin{equation}\label{eq:W-transformation-matrix}
        \left(\mathcal{W}_s\right)_{i,j} = w_{i,j} = P_{j-1}(c_i), \quad i=1,\ldots,k, \qquad j=1,\ldots,s.
    \end{equation}
In~\eqref{eq:what-we-have-to-compute} we, therefore, recognize a method that can be expressed with the Butcher Tableau~\eqref{eq:butcher_tableau}
\begin{equation}\label{eq:rank-deficient-tableau}
    \mathbf{c}^\top = [c_1,\ldots,c_k], \quad A = \mathcal{I}_s \mathcal{W}_s^\top \mathcal{B} \in \mathbb{R}^{k \times k}, \quad %
    \mathbf{b}^\top = [b_1,\ldots,b_k], %
\end{equation}
where
\begin{equation}\label{eq:ps_is_b_matrices}
\begin{split}
\mathcal{I}_s = & \; \begin{bmatrix}
    \int_{0}^{c_1} P_0(x)\,\mathrm{d}x & \cdots & \int_{0}^{c_1} P_{s-1}(x)\,\mathrm{d}x\\
    \vdots & \ddots & \vdots \\
    \int_{0}^{c_k} P_0(x)\,\mathrm{d}x & \cdots & \int_{0}^{c_k} P_{s-1}(x)\,\mathrm{d}x
\end{bmatrix}_{k \times s}, %
\quad \mathcal{B} = \begin{bmatrix}
    b_1 & \\
    & \ddots & \\
    & & b_k
\end{bmatrix},
\end{split}
\end{equation}
and $\mathcal{W}_s$ is the transformation matrix in~\eqref{eq:W-transformation-matrix}; which is a RK
formulation of the form~\eqref{eq:steps_Runge_Kutta}-\eqref{eq:stages_Runge_Kutta_v2} of such method. 
We immediately observe from~\eqref{eq:rank-deficient-tableau} that the matrix $A$ of the Butcher Tableau is rank deficient, that is, by construction it has maximum rank $s \leq k$. Hence, we can restrict the computation in~\eqref{eq:what-we-have-to-compute} to only $s$ equations and recover the remaining $k-s$ stages as a linear combination of the others. To avoid fixing arbitrarily a subset of the stages, we can recast the problem as computing the $\{ \widehat{\boldsymbol{\phi}}_\ell \}_{\ell=0}^{s-1}$ in~\eqref{eq:relation-for-the-phi} which are in the right number and do not require us to make a choice. For all time-steps, the resulting set of equations is then given by
\begin{equation}\label{eq:the-formulation-we-use}
\begin{split}
    \widehat{\boldsymbol{\phi}}_j = &\; \sum_{i=1}^{k} b_i P_j(c_i) f\left( \mathbf{y}^{(n)} + h \sum_{\ell = 0}^{s-1} \widehat{\boldsymbol{\phi}}_\ell \int_{0}^{c_i} P_\ell(x)\,\mathrm{d}x \right), \; j=0,\ldots,s-1\\
    \mathbf{y}^{(n+1)} = &\; \mathbf{y}^{(n)} +  h\, \widehat{\boldsymbol{\phi}}_0.
\end{split}
\end{equation}
From this, we have obtained the set of exactly $s$ nonlinear equations whose solution is needed to advance from one time-step to the following. 

To construct and implement the method described above, it is convenient to express its coefficient matrices in a compact form, analogous to the formulation used for RK collocation methods; see, e.g., \cite[\S IV.5]{BibleVolII}. This can be done thanks to the following result~\cite{BrugnanoBook}.
 \begin{proposition}[{Matrix Structure, \cite[Theorem 3.7]{BrugnanoBook}}]\label{pro:matrix_structure}
Consider the $\operatorname{HBVM}(k,s)$ scheme, $k \geq s$, the matrix $\mathcal{W}_s$, $\mathcal{I}_s$ and $\mathcal{B}$ in~\eqref{eq:W-transformation-matrix}--\eqref{eq:ps_is_b_matrices} are such that
\[
\mathcal{I}_s = \mathcal{W}_{s+1} \hat{X}_s, \quad \hat{X}_s = \begin{bmatrix}
    \nicefrac{1}{2} & - \xi_1 \\
    \xi_1 & & \ddots \\
    & \ddots & & -\xi_{s-1} \\
    & & \xi_{s-1} \\
    \hline
    & & & \xi_{s} 
\end{bmatrix} = \begin{bmatrix}
        & X_s & \\
        \hline
        0 & \cdots & \xi_s
    \end{bmatrix}, \; \xi_i = \frac{1}{2\sqrt{4 i^2 - 1}},
\]
hence $X_s = \mathcal{W}_s^\top \mathcal{B} \mathcal{I}_s$.
\end{proposition}

\section{The linear case}\label{sec:thelinearcase}
To discuss the procedure for solving~\eqref{eq:the-formulation-we-use} we start from the simplifying assumption on the linearity of $f$, i.e.,
\begin{equation}\label{eq:linearity_assumption}
    f(\mathbf{y}) = \mathcal{G} \mathbf{y}, \quad \mathcal{G} \in \mathbb{R}^{2m \times 2m}.
\end{equation}
Under this hypothesis, we can restate~\eqref{eq:the-formulation-we-use} as a linear matrix equation. We define the matrix 
\[
\Phi = \left[ \widehat{\boldsymbol{\phi}}_0 | \ldots | \widehat{\boldsymbol{\phi}}_{s-1} \right] \in \mathbb{R}^{2m \times s},
\]
collecting all the stages in its columns, and then state the matrix equation as
\begin{equation}\label{eq:matrix-eq-formulation}
\Phi - h \mathcal{G} \Phi X_s^\top = \mathcal{G} \mathbf{y}^{(n)} \mathbf{b}^\top \mathcal{W}_s, 
\end{equation}
where we have used that $X_s^\top = \mathcal{I}_s^\top \mathcal{B} \mathcal{W}_s$, see Proposition~\ref{pro:matrix_structure}. The previous can be rewritten as the linear system
\begin{equation}\label{eq:linear-system-formulation}
\left[ I - h X_s \otimes \mathcal{G}  \right] \boldsymbol{\varphi} = \mathcal{W}_s^\top \mathbf{b} \otimes \mathcal{G} \mathbf{y}^{(n)},
\end{equation}
where $\boldsymbol{\varphi} = \operatorname{vec}(\Phi)$ is obtained by stacking the columns of $\Phi$, $I$ is the $2sm \times 2sm$ identity matrix, and $\otimes$ represents the Kronecker product. To compute the reduced formulation stages, we must solve the matrix equation~\eqref{eq:matrix-eq-formulation}. One approach is to solve the corresponding larger linear system in~\eqref{eq:linear-system-formulation}, nevertheless it is better to exploit the structure of the matrix equation~\eqref{eq:matrix-eq-formulation}. Observe that the right-hand side is a rank-1 matrix, indeed
\[
\mathcal{G} \mathbf{y}^{(n)} \mathbf{b}^\top \mathcal{W}_s = \mathbf{f}_1 \mathbf{f}_2^\top, \text{ with } \mathbf{f}_1 = \mathcal{G} \mathbf{y}^{(n)} \in \mathbb{R}^{2m}, \quad \mathbf{f}_2 = \mathcal{W}_s^\top \mathbf{b} \in \mathbb{R}^s. 
\]
To solve~\eqref{eq:matrix-eq-formulation} efficiently, we adopt the one-sided Krylov projection method~\cite{SimonciniReview}, which constructs a low-dimensional approximation of the solution by iteratively projecting onto a Krylov subspace. Specifically, we can further rewrite~\eqref{eq:matrix-eq-formulation} as
\begin{equation}\label{eq:sylvester}
    Z E + E R = \mathbf{u} \mathbf{v}^\top, \quad \text{where} \quad 
    \begin{cases}
        Z = -h^{-1} \mathcal{G}^{-1},\\
        R = X_s^{\top},\\
        \mathbf{u} = -Z \mathbf{f}_1,\\
        \mathbf{v} = \mathbf{f}_2,
    \end{cases}
\end{equation}
and consider the space generated by the large matrix $Z$ and the vector $\mathbf{u}$. This reduces the problem to solving a Sylvester equation of the same form but with reduced dimension. Once the reduced system is solved, the approximate solution to the original problem is reconstructed by lifting back to the original space.
{\begin{remark}\label{rmk:sparse-dense}
The choice of Krylov subspace methods for the solution of equation in~\eqref{eq:sylvester} is not the only available one. Equation~\eqref{eq:sylvester} is an example of \emph{sparse-dense} Sylvester equation~\cite{SparseDense} in which we have the two coefficients being a smaller dense matrix $R$ and a large sparse matrix $Z$. To solve such a problem one could consider diagonalizing $R$ and then solving a sequence of shifted linear systems~\cite{MR1917499}. The stability of such procedure depends on the conditioning of the eigenvectors matrix. Observe that by construction $X_s$ is non normal, and it is known that the conditioning of its eigenvectors grows exponentially with $s$; see, e.g.,~\cite{durastante2025stageparallelimplicitrungekuttamethods}. This implies that for larger values of $s$ the procedure could incur in numerical instabilities. To accommodate for this, the approach in~\cite{durastante2025stageparallelimplicitrungekuttamethods} was to perturb the $R$ matrix and build a two-step diagonal-solution/low-rank correction procedure. In the present case, since the right-hand side is already low-rank, one can skip the perturbation step, and proceed through with the low-rank solution approach.     
\end{remark}}
The {solution} procedure is presented in its simplest form in Algorithm~\ref{alg:oneside_sylvester}, where a polynomial Krylov subspace is used to project the Sylvester equation onto a lower-dimensional space.
\begin{algorithm}[htbp]
\KwData{Matrices $Z \in \mathbb{R}^{N \times N}$, $R \in \mathbb{R}^{s \times s}$, vectors $\mathbf{u} \in \mathbb{R}^{N}$, $\mathbf{v} \in \mathbb{R}^s$, tolerance $\epsilon$, maximum iterations $k_{\max}$}
\KwResult{Approximate solution $E$ to $Z E + E R = \mathbf{u} \mathbf{v}^\top$}

\BlankLine
Initialize $V_1 = \frac{\mathbf{u}}{\|\mathbf{u}\|}$, $\beta = \|\mathbf{u}\|$\;
\For{$j = 1, 2, \dots, k_{\max}-1$}{
    \tcc{Krylov subspace generation}
    Compute $\mathbf{w} = Z \mathbf{v}_j$\;
    \For{$i = 1, 2, \dots, j$}{
        Compute $h_{i,j} = \mathbf{v}_i^\top \mathbf{w}$\;
        Orthogonalize: $\mathbf{w} = \mathbf{w} - h_{i,j} \mathbf{v}_i$\;
    }
    Compute $h_{j+1,j} = \|\mathbf{w}\|$\;
    \If{$h_{j+1,j} = 0$}{
        \textbf{break}\;
    }
    Normalize: $\mathbf{v}_{j+1} = \frac{\mathbf{w}}{h_{j+1,j}}$\;
    Update $V_{j+1} = [\mathbf{v}_1| \cdots | \mathbf{v}_{j+1}]$\;
Form the matrix $H_j = V_j^\top Z V_j$\;
\tcc{Projected problem solution}
Solve the projected Sylvester equation $H_j Y_j + Y_j R = \beta e_1 \mathbf{v}^\top$\;
Compute the residual: $\rho_j = |h_{j+1,j}| \| \mathbf{e}_j^\top Y_j \|_F$\;
\If{$ \rho_j < \epsilon \times \| \mathbf{u} \|_F \| \mathbf{v} \|_F$}{
    \textbf{break}\;
}
}
Compute the approximate solution $E_j = V_j Y_j$\;
\caption{One-sided Krylov projection method for the matrix equation $Z E + E R = \mathbf{u} \mathbf{v}^\top$ with Arnoldi iteration}\label{alg:oneside_sylvester}
\end{algorithm}
Specifically, the projection is performed onto the space $\mathcal{K}_i(Z,\mathbf{u}) = \operatorname{span}\left\{\mathbf{u}, Z\mathbf{u}, Z^2\mathbf{u}, \ldots, Z^{i-1}\mathbf{u} \right\},$ where $i$ denotes the dimension of the Krylov subspace. While the use of a polynomial Krylov subspace provides a straightforward way for approximating the solution of the Sylvester equation, its convergence can be slow, particularly when the spectrum of \( Z \) is wide or ill-conditioned. To address this limitation and accelerate convergence, one may employ an \emph{extended Krylov subspace method}~\cite{MR2318706}, which augments the approximation space by including both powers and inverses of \( Z \), and is defined as
\begin{equation}\label{eq:extended_krylov}
\mathcal{EK}_i(Z,\mathbf{u}) = \operatorname{span}\left\{ \mathbf{u}, Z^{-1}\mathbf{u}, Z\mathbf{u}, Z^{-2}\mathbf{u}, \ldots, Z^{-(i-1)}\mathbf{u}, Z^{i-1}\mathbf{u} \right\},    
\end{equation}
thereby capturing a broader range of spectral components. This enhanced subspace typically leads to significantly faster convergence{. The gain in convergence is usually paid by the need for the solution of linear systems with $Z$. In our case, this computation is reduced to the cost of a matrix-vector product with the matrix \( \mathcal{G} \), and the true cost appears instead in the computation of the matrix-vector products with \( Z = -h^{-1} \mathcal{G}^{-1} \).}
Finally, one can also employ more general rational Krylov subspace methods for~\eqref{eq:matrix-eq-formulation}, by constructing approximation spaces of the form $\mathrm{span}\{(Z-\xi_1 I)^{-1}\mathbf{u},\ldots,(Z-\xi_k I)^{-1}\mathbf{u}\}$ with selected shifts (poles)~$\xi_i$. The convergence of such rational-Krylov methods critically depends on the choice of these poles, and various algorithms have been proposed to select or adapt them. For example,~\cite{CasulliRobol2023} develops an adaptive block rational Krylov solver for Sylvester equations that keeps the last pole at infinity and updates the others based on convergence estimates, while~\cite{BergermannStoll2023} discusses adaptive pole selection strategies (and error estimates) in the context of exponential integrators. Techniques like the RKFIT algorithm~\cite{Guettel2017} can also be used to compute near-optimal poles by fitting a rational function to the matrix exponential over a spectral interval. When poles are chosen to target the spectrum of $Z$ (for instance via minimax rational approximation), rational-Krylov methods typically converge far faster than polynomial Krylov methods (see, e.g., \cite{Beckermann2011,BergermannStoll2023}). However, these improvements come at the cost of solving a sequence of shifted linear systems $(Z - \xi_i I)\mathbf{w} = \mathbf{u}$, whose numerical difficulty depends strongly on the magnitude and distribution of the poles and may require robust preconditioners or factorization reuse strategies to be computationally efficient. {Alternative approaches to the Krylov settings are represented by alternating-direction implicit methods~\cite{MR2557293} which still require solving shifted linear systems with prescribed poles, which can be modified to exploit low-rank structure in the right-hand side of the equation~\cite{MR3197189}, and the direct method exploiting the sparse-dense formulation~\cite{SparseDense} of~\eqref{eq:sylvester} already mentioned in Remark~\ref{rmk:sparse-dense}. For the latter, this would increase both the memory requirements, and the number of factorizations to be computed with respect to the effort needed to expand the Extended Krylov subspace~\eqref{eq:extended_krylov}.}

\section{The nonlinear case}\label{sec:thenonlinearcase}

To apply the procedure to the nonlinear case, we start from the formulation in~\eqref{eq:the-formulation-we-use} and apply a Newton-like linearization. We define the nonlinear system for $j=0,\ldots,s-1:$
\begin{equation}\label{eq:nonlinear-system}
\mathcal{F}_j(\widehat{\boldsymbol{\phi}}_0, \ldots, \widehat{\boldsymbol{\phi}}_{s-1}) := \widehat{\boldsymbol{\phi}}_j - \sum_{i=1}^{k} b_i P_j(c_i)\, f\left( \mathbf{y}^{(n)} + h \sum_{\ell = 0}^{s-1} \widehat{\boldsymbol{\phi}}_\ell \int_{0}^{c_i} P_\ell(x)\,\mathrm{d}x \right),
\end{equation}
and aim to solve for the matrix $\Phi \in \mathbb{R}^{2m \times s}$ collecting all the \(\widehat{\boldsymbol{\phi}}_0, \ldots, \widehat{\boldsymbol{\phi}}_{s-1}\) on the columns. We rewrite~\eqref{eq:nonlinear-system} in matrix format as
\begin{equation}\label{eq:f=0_in_matrix_form}
\mathcal{F}(\Phi) = \Phi  - \left[ \left. f\left( \mathbf{y}^{(n)} + h \Phi \mathcal{I}_s^\top \mathbf{e}_1  \right)  \right| \cdots \left|  f\left( \mathbf{y}^{(n)} + h \Phi \mathcal{I}_s^\top \mathbf{e}_{k}  \right) \right. \right]_{2m \times k} \mathcal{B} \mathcal{W}_s = 0,    
\end{equation}
from which ones recovers the matrix equation in~\eqref{eq:matrix-eq-formulation} under the linearity assumption~\eqref{eq:linearity_assumption} for $f$. If we apply a \emph{simplified Newton} approach~\cite{MR1790038} for solving $\mathcal{F}(\Phi) = 0$ as in~\cite[\S 4.2]{BrugnanoBook} we obtain\footnote{With respect to the original source, we have adapted the notation to avoid evaluation of vector function on block-vectors by means of the Kronecker product; see the discussion in~\cite[\S 4, Eq.~(4.2)]{BrugnanoBook}.} an iteration for the matrix increment sequence $\{\Delta^{[p]} \in \mathbb{R}^{2m \times s}\}_p $ of the form
\begin{equation}\label{eq:matrix_form_of_simplified_Newton}
\begin{cases}
\Delta^{[p]} - h J_f \Delta^{[p]} \mathcal{I}_s^\top \mathcal{B} \mathcal{W}_s =  \mathcal{F} (\Phi^{[p]}),  \\
\Phi^{[p+1]} = \Phi^{[p]} - \Delta^{[p]}.
\end{cases}
\end{equation}
Here, the \emph{simplification} has come by having substituted to the Jacobians $\{  J_{f,j} = \left.\nicefrac{\partial f}{\partial \mathbf{y}}\right|_{\mathbf{y} = \widehat{\boldsymbol{\phi}}_{j}} \}_{j=0}^{s-1}$---which would be needed for expressing the full Jacobian $J_{\mathcal{F}}$ of~\eqref{eq:f=0_in_matrix_form} in a complete Newton iteration---the approximated Jacobian $J_{f} = \left.\nicefrac{\partial f}{\partial \mathbf{y}}\right|_{\mathbf{y} = \mathbf{y}^{(n)}}$. 

Observe that by construction $\mathcal{F}(\Phi^{[p]}) \in \mathbb{R}^{2m\times s}$ with $s \ll m$, hence it is a \emph{low-rank} matrix. To recover $\Delta^{[p]}$ we start by decomposing 
\[
\mathcal{F}(\Phi^{[p]}) = F_1 F_2^\top, \qquad F_1 \in \mathbb{R}^{2m \times r}, \; F_2 \in \mathbb{R}^{s \times r}, \qquad r = \operatorname{rank}(\mathcal{F}(\Phi^{[p]})) \leq s,
\]
and identifying $\mathcal{I}_s^\top \mathcal{B} \mathcal{W}_s = X_s^\top$ (cf. Proposition~\ref{pro:matrix_structure}). We can therefore solve the Sylvester matrix equation
\[
h^{-1} J_f^{-1} \Delta^{[p]} - \Delta^{[p]} X_s^\top = \tilde{F}_1 F_2^\top, \qquad \tilde{F}_1 = h^{-1} J_f^{-1} F_1,
\]
and finally compute the update in~\eqref{eq:matrix_form_of_simplified_Newton}.

To enhance the convergence of the simplified Newton method, which is generally less robust than the one for the full Newton---possibly even requiring restriction on the time-step, it is often beneficial to first perform a few fixed-point iterations to produce a more accurate initial guess for the Newton solver. Specifically, starting from \(\Phi^{(0)} = 0\), we compute
\begin{equation}\label{eq:fixed_point_initialization}
\Phi^{(\ell)} = 
\left[\, f\Big( \mathbf{y}^{(n)} + h \Phi^{(\ell-1)} \mathcal{I}_s^\top \mathbf{e}_1 \Big) \;\big|\; \cdots \;\big|\; f\Big( \mathbf{y}^{(n)} + h \Phi^{(\ell-1)} \mathcal{I}_s^\top \mathbf{e}_k \Big) \,\right]_{2m \times k} \; \mathcal{B} \mathcal{W}_s,
\end{equation}
for \(\ell = 1, \dots, \ell_{\rm max}\).  
The output \(\Phi^{(\ell_{\rm max})}\) is then used as the starting point for the simplified Newton iteration. This procedure effectively \emph{warms up} the iteration, bringing the initial guess closer to the true solution and reducing the number of subsequent Newton steps required for convergence.

\subsection{Preconditioner for the complete Newton--Krylov case}\label{sec:preconditioner}

In the complete (full) Newton case we do not replace the stage Jacobians by a single frozen Jacobian. Consequently, the Kronecker product structure that permitted us to write the solution of the linear system in the matrix equation form~\eqref{eq:matrix_form_of_simplified_Newton} is lost, and one must work with the full block Jacobian of the nonlinear map \(\mathcal{F}:\mathbb{R}^{2m\times s}\to\mathbb{R}^{2m\times s}\) introduced in~\eqref{eq:f=0_in_matrix_form}.  

To write this Jacobian in a convenient form, let
\[
J_{f,i} \;=\; \left.\frac{\partial f}{\partial \mathbf y}\right|_{\mathbf y=\mathbf{y}^{(n)} + h\sum_{\ell=0}^{s-1}\widehat{\boldsymbol\phi}_\ell \int_0^{c_i} P_\ell(x)\,\mathrm dx}\in\mathbb R^{2m\times 2m},
\qquad i=1,\dots,k,
\]
be the stage-dependent Jacobians, and define
\[
S_{i,\ell}=\int_0^{c_i} P_\ell(x)\,\mathrm{d}x,\, i=1,\dots,k,\; \ell=0,\dots,s-1.
\]
By differentiating~\eqref{eq:nonlinear-system} with respect to the unknown columns \(\widehat{\boldsymbol\phi}_\ell\), one obtains that the full Newton Jacobian \(J_{\mathcal F}\) is an \(s\times s\) block matrix with \(2m\times 2m\) blocks:
\begin{equation}\label{eq:full-block-Jacobian}
\big(J_{\mathcal F}\big)_{j,\ell}
\;=\;
\frac{\partial \mathcal F_j}{\partial \widehat{\boldsymbol\phi}_\ell}
\;=\;
\delta_{j\ell}\,I_{2m}
\;-\; h\sum_{i=1}^{k} b_i\,P_j(c_i)\; J_{f,i}\; S_{i,\ell},
\qquad j,\ell=0,\dots,s-1,
\end{equation}
where \(\delta_{j\ell}\) is the Kronecker delta. Equivalently, if we vectorize the unknowns column-wise, \( \operatorname{vec}(\Phi)\in\mathbb R^{2ms}\), then \(J_{\mathcal F}\) is a \(2ms\times 2ms\) matrix built from the blocks in~\eqref{eq:full-block-Jacobian}.  
Because \(J_{f,i}\) depends on the stage \(i\), the double summation in~\eqref{eq:full-block-Jacobian} prevents a representation of \(J_{\mathcal F}\) as a simple Kronecker product of a small \(s\times s\) matrix with a single \(2m\times 2m\) matrix; this is precisely the loss of Kronecker structure.
Furthermore, assembling and storing the full Jacobian~\eqref{eq:full-block-Jacobian} comes at increased cost per iteration of the Newton method. A viable alternative is represented by the usage of a finite difference approximation of the action of the Jacobian $J_{\mathcal F}$ against a vector $\mathbf{w}$, which permits us to implement a Newton--Krylov procedure, i.e., a Newton method in which the linear systems with the Jacobian are solved via a Krylov method needing only access to the matrix-vector product routine. Hence, to avoid forming the full Jacobian matrix of the nonlinear system~\eqref{eq:f=0_in_matrix_form}, we approximate the product of the Jacobian of \(\mathcal{F}:\mathbb{R}^{2m\times s}\to\mathbb{R}^{2m\times s}\) with a given perturbation \(\mathbf{W}\in\mathbb{R}^{2m\times s}\) by means of the following finite-difference directional derivative.  
For a given iterate \(\Phi\in\mathbb{R}^{2m\times s}\), we compute
\begin{equation}\label{eq:fd-dirder-matrix}
J_{\mathcal{F}}(\Phi)\,\mathbf{W} 
\;\approx\;
\frac{\mathcal{F}(\Phi+\varepsilon\,\mathbf{W}) - \mathcal{F}(\Phi)}{\varepsilon},
\end{equation}
where \(\varepsilon>0\) is a suitably scaled perturbation parameter.  
In practice, \(\varepsilon\) is chosen as
\begin{equation}\label{eq:eps-scaling}
\varepsilon \;=\; \frac{10^{-7}}{\|\mathbf{W}\|_F}\,
\max\!\big(1,\,\|\Phi\|_F\big),
\end{equation}
so that the step size is adaptively scaled both with the Frobenius norm of the current iterate and with the magnitude of the perturbation direction.  
This scaling ensures numerical robustness and prevents underflow or overflow in finite precision.  
If \(\|\mathbf{W}\|_F = 0\), the routine simply returns the zero matrix. The resulting approximation~\eqref{eq:fd-dirder-matrix} provides a first-order accurate estimate of the Jacobian action. Observe also that each new Jacobian evaluation comes at the cost of a single evaluation of $\mathcal{F}$ on $\Phi+\varepsilon\,\mathbf{W}$, since we have stored $\mathcal{F}(\Phi)$, which is also the right-hand side of the system, and the current value of the approximation.

Since, in the general case, the Jacobian matrix associated with the nonlinear system~\eqref{eq:full-block-Jacobian} is \emph{non-symmetric} and lacks any special structure, an appropriate choice for the iterative linear solver is required. Among Krylov subspace methods~\cite{Saad}, the most natural candidates for handling non-symmetric, potentially ill-conditioned systems are the \emph{Generalized Minimal Residual method} (GMRES)~\cite{MR848568} and its restarted variant GMRES$(m)$, as well as the \emph{Bi-Conjugate Gradient Stabilized method} (BiCGStab)~\cite{MR1283335}.  
These algorithms are well suited for large-scale, sparse, or matrix-free contexts such as ours, where the Jacobian is never assembled explicitly and only its action on test matrices or vectors is available via the finite-difference approximation~\eqref{eq:fd-dirder-matrix}.

In practice, the convergence of such iterative solvers is strongly influenced by the conditioning of the linear system, which deteriorates as the time-step \(h\) increases or as the stiffness of the underlying differential problem grows.  To alleviate this, we need to introduce a preconditioner approximating the inverse of the Jacobian \(J_{\mathcal{F}}(\Phi)\). The preconditioner can be constructed based on the simplified Newton structure and applied at each iteration to accelerate convergence; specifically, if we replace every stage Jacobian \(J_{f,i}\) by the Jacobian at the previous Newton or step iteration, \(J_f:=\left.\partial f/\partial\mathbf y\right|_{\mathbf y=\mathbf y^{(n)}}\). Then
  \[
  \big(J_{\mathcal F}^{\mathrm{(approx)}}\big)_{j,\ell}
  =\delta_{j\ell}I_{2m}-h\sum_{i=1}^k b_i P_j(c_i)\,J_f\,S_{i,\ell}
  =\delta_{j\ell}I_{2m}-h\big(X_s\big)_{j,\ell}\,J_f,
  \]
  where \(X_s=\mathcal W_s^\top\mathcal B\mathcal I_s\) (cf.\ Proposition~\ref{pro:matrix_structure}). In vectorized form this yields the matrix equation representation
  \begin{equation}\label{eq:matrix_equation_for_Newton}
    W - h J_f W X_s^\top =  R \equiv R_1 R_2^\top, \qquad R_1 \in \mathbb{R}^{2m \times r},\; \quad R_2 \in \mathbb{R}^{s \times r}, \qquad r = \operatorname{rank}(R),    
  \end{equation}
  where $W \in \mathbb{R}^{2m \times s}$ represents the preconditioned residual, and $R \in \mathbb{R}^{2m \times s}$ is the residual at the current step. The equation~\eqref{eq:matrix_equation_for_Newton} can be solved more cheaply than the full Jacobian by means of an inner one-sided projection method as in Algorithm~\ref{alg:oneside_sylvester} and therefore serves as an effective preconditioner. This is the natural extension of the simplified Newton Jacobian used previously. We stress that using Algorithm~\ref{alg:oneside_sylvester} introduces a variability in the preconditioner that violates the fixed-preconditioner assumption of standard GMRES, and therefore necessitates the use of the \emph{Flexible GMRES} (FGMRES) algorithm~\cite{MR1204241}, which is specifically designed to accommodate changing or nonlinear preconditioners, which in our setting provides a robust and efficient matrix-free framework for solving the linearized Newton systems.

  {Algorithm~\ref{alg:complete_Newton_Krylov} makes this construction explicit. At each Newton step, the Jacobian action is approximated through the directional finite difference~\eqref{eq:fd-dirder-matrix}, while the linearized system is solved inexactly by FGMRES. Each application of the right preconditioner requires solving~\eqref{eq:matrix_equation_for_Newton} with the frozen Jacobian $J_f$, so the preconditioner is refreshed together with the current Newton iterate and remains consistent with the local linearization.}

\begin{algorithm}[htbp]
\KwData{State $\mathbf{y}^{(n)}$, initial step size $h$, initial guess $\Phi^{[0]}$, Newton tolerance $\tau_{\rm newton}$, maximum Newton iterations $p_{\max}$, initial inner FGMRES tolerance $\eta_0$, upper bound $\eta_{\max}$}
\KwResult{Accepted stage matrix $\Phi$ and updated step size $h$}
\Repeat{a Newton solve converges within $p_{\max}$ iterations}{
    Compute a warm start $\Phi^{[0]}$ by the fixed-point iteration~\eqref{eq:fixed_point_initialization}\;
    Set $p=0$ and $f_{\mathrm{nrmo}}=\infty$\;
    \While{$p < p_{\max}$}{
        Compute $R^{[p]} = -\mathcal{F}(\Phi^{[p]})$ and $f_{\mathrm{nrm}} = \|R^{[p]}\|_F$\;
        \If{$f_{\mathrm{nrm}} \le \tau_{\rm newton}$}{
            accept the step and \textbf{break}\;
        }
        Define the Jacobian action $J_{\mathcal F}(\Phi^{[p]})[\cdot]$ by~\eqref{eq:fd-dirder-matrix}\;
        Approximately solve
        \[
        J_{\mathcal F}(\Phi^{[p]})\,\Delta^{[p]} = R^{[p]}
        \]
        by right-preconditioned FGMRES with stopping tolerance $\eta_p$, where each preconditioner application computes $W$ from~\eqref{eq:matrix_equation_for_Newton} by Algorithm~\ref{alg:oneside_sylvester}\;
        Update the Newton iterate $\Phi^{[p+1]} = \Phi^{[p]} + \Delta^{[p]}$\;
        Update the inner FGMRES tolerance $\eta_{p+1}$ using the residual ratio strategy described in Section~\ref{sec:deciding_tolerances}\;
        Set $f_{\mathrm{nrmo}} = f_{\mathrm{nrm}}$ and $p = p+1$\;
    }
    \If{$p = p_{\max}$ and $\|\mathcal{F}(\Phi^{[p]})\|_F > \tau_{\rm newton}$}{
        reject the step, set $h = h/2$, and restart the Newton solve with the reduced step size\;
    }
}
\caption{{Inexact Newton--Krylov iteration for the HBVM stage equations with FGMRES and matrix-equation preconditioning}}
\label{alg:complete_Newton_Krylov}
\end{algorithm}

\subsection{Setting the tolerances for the Newtok--Krylov method}\label{sec:deciding_tolerances}
    In the inexact Newton–Krylov framework, as discussed in~\cite[\S 6.2.1]{MR1344684}, a crucial design decision concerns how tightly the inner linear solver---FGMRES in our case---must converge before proceeding with the next outer Newton step.  If the inner residual satisfies  
\[
\| \mathcal{F}(\Phi^{[p]}) + J_{\mathcal F}(\Phi^{[p]})\,\Delta^{[p]} \|_F \;\le\; \eta_p \,\| \mathcal{F}(\Phi^{[p]})\|_F,
\]
for some forcing parameter \(\eta_p \in [0,\eta_{\max})\), then the overall inexact Newton iteration converges under the usual assumptions, i.e., that the nonlinear mapping \(\mathcal{F}\) is continuously differentiable in a neighbourhood of the solution \(\Phi^{*}\), that its Jacobian \(J_{\mathcal{F}}(\Phi^{*})\) is nonsingular, and that \(J_{\mathcal{F}}\) satisfies a Lipschitz, or H\"{o}lder, continuity condition in that neighbourhood.  
Under these hypotheses, the inexact Newton iteration retains global convergence properties and recovers the classical local quadratic (or superlinear) convergence rate as the forcing parameter \(\eta_p \to 0\). The parameter \(\eta_p\) should be chosen adaptively: large enough when \(\|\mathcal{F}(\Phi^{[p]})\|_F\) is still large to avoid \emph{oversolving} the linear system, yet decreasing as the Newton iteration progresses and the nonlinear residual becomes small.  By coupling the inner solver tolerance \(\tau_{\rm lin}\approx\eta_p\,\|\mathcal{F}(\Phi^{[p]})\|_F\) with the outer Newton tolerance \(\tau_{\rm newton}\), one achieves a balance between efficiency and robustness: the linear solver is not forced to extreme accuracy when heavy correction is still required, yet is tightened in the final stages to preserve quadratic or super-linear convergence.

In the actual implementation, the convergence of the inner FGMRES solver is dynamically controlled based on the progress of the outer Newton iteration, rather than using a fixed forcing term.  Let \(\mathbf{f}_0 = \operatorname{vec}(\mathcal{F}(\Phi^{[p]}))\) and \(f_{\text{nrm}} = \|\mathbf{f}_0\|/\sqrt{n}\) denote the normalized residual at the current Newton step. After each outer iteration, the ratio of successive residual norms, 
\(\mathrm{rat} = f_{\text{nrm}} / f_{\text{nrmo}}\), is computed, where \(f_{\text{nrmo}}\) is the residual norm from the previous Newton step. This ratio is then used to update the effective linear-solver tolerance \(\eta_p\) for the next step according to
\[
\eta_p = \min\Big\{ \eta_{\text{max}}, \max\big( \gamma \, \eta_{\rm old}^2, \, \gamma \, \mathrm{rat}^2, \, 0.5 \, \mathrm{stop\_tol}/f_{\text{nrm}} \big) \Big\},
\]
where \(\eta_{\rm old}\) is the previous tolerance, \(\gamma\) is a user-specified parameter controlling the aggressiveness of the update, \(\eta_{\text{max}}\) is an upper bound on the forcing term, and the stop tolerance is computed as $\mathrm{stop\_tol}= \varepsilon_{\text{abs}} + \varepsilon_{\text{rel}} \times f_{\text{nrm}}$ for $(\varepsilon_{\text{abs}},\varepsilon_{\text{rel}})$ the absolute e relative tolerances given by the user. This adaptive mechanism ensures that the inner GMRES solve is not oversolved when the Newton residual is large, but the tolerance is tightened as the iteration progresses and the solution is approached; we refere again to~\cite[\S 6.2.1]{MR1344684} for a complete discussion on the construction of the Newton--Krylov method. 

As for the simplified Newton case, we can use~\eqref{eq:fixed_point_initialization} to \emph{warm start} the complete Newton iteration. To further enhance the robustness of the Newton--Krylov solver, particularly for challenging nonlinear problems, we employ a simple time-step adaptivity strategy that ensures convergence. If the Newton iteration fails to converge within a prescribed maximum number of iterations, the time-step is halved and the computation for that time-step is restarted. Conversely, when the Newton method converges successfully for several consecutive steps, the time step is gradually increased to improve computational efficiency. Specifically, if the Newton iteration converges successfully for more than 4 consecutive time steps, the step size $h$ is doubled, up to a maximum value specified by the user. This adaptive mechanism provides a balance between robustness---by reducing the step size when convergence difficulties arise---and efficiency---by increasing the step size when the solution evolves smoothly.

\section{Numerical examples}
\label{sec:numerical_examples}

For the numerical experiments, we consider a semi-discretization in space of a partial differential equation with Hamiltonian structure. We begin by examining the linear case to demonstrate the performance of the matrix equation solver introduced in Section \ref{sec:thelinearcase}. For the nonlinear case, our attention centers on the preconditioner described in Section~\ref{sec:preconditioner}. We do not include a separate, detailed analysis of the simplified Newton method, as its implementation effectively reduces to repeated applications of the linear case; its practical impact is therefore captured through the performance of the preconditioner itself (indeed, every application of the preconditioner coincides with a simplified Newton step).

All numerical experiments in this section were performed on a single node of the Toeplitz cluster, located at the Green Data Center of the University of Pisa. The node is equipped with an Intel\textsuperscript{\textregistered} Xeon\textsuperscript{\textregistered} CPU E5-2643 v4 at 3.40\si{\giga\hertz} with 2 threads per core, 6 cores per socket and 2 sockets, with 128\si{\giga\byte} of RAM memory. We run the experiments using MATLAB \texttt{v9.10.0.1602886} (R2021a). For the construction of the HBVM($k$,$s$) matrices described in Proposition~\ref{pro:matrix_structure} we have used the routine \lstinline[language=Matlab]{RKform(k,s)} from~\cite{MR2833606}, while for the implementation of the one-sided Krylov projection method for the solution of the matrix equations~\eqref{eq:sylvester} and~\eqref{eq:matrix_equation_for_Newton} we have modified the original \lstinline[language=Matlab]{kpik} solver from~\cite{MR2318706}.  

\subsection{A semi-discretized linear wave equation}\label{sec:semidiscrete_linear}

We apply the strategy from Section~\ref{sec:thelinearcase} to the one-dimensional linear wave equation
\begin{equation*}
    u_{tt}(x,t) = u_{xx}(x,t), \qquad x\in[0,L],\; t\ge 0,
\end{equation*}
where $u(x,t)$ denotes the displacement field. Introduce a uniform spatial grid $x_j=j {\delta x}$ for $j=0,1,\dots,N$ with ${\delta x}=\nicefrac{L}{N}$, and denote
$u_j(t)\approx u(x_j,t)$.  Approximating the second derivative by central finite differences gives, for interior nodes $j=1,\dots,N-1$,
\begin{equation*}
    u_{tt}(x_j,t) \approx \,\frac{u_{j+1}(t)-2u_j(t)+u_{j-1}(t)}{{\delta x}^{2}}.
\end{equation*}
Collecting the interior degrees of freedom into the vectors
\[
\mathbf{u}(t)=(u_1(t),\dots,u_{N-1}(t))^\top,\qquad 
\mathbf{p}(t)=\dot{\mathbf{u}}(t),
\]
we form the combined state vector
\begin{equation}\label{eq:state-vector}
        \mathbf{y}(t)=\begin{pmatrix}\mathbf{u}(t)\\[2pt]\mathbf{p}(t)\end{pmatrix}.
\end{equation}
Let $L_N\in\mathbb{R}^{(N-1)\times(N-1)}$ denote the standard discrete Laplacian (with homogeneous Dirichlet boundary conditions at $x=0$ and $x=L$),
\[
    L_N=\frac{1}{{\delta x}^{2}}
    \begin{bmatrix}
    -2 & 1  &        &        &   \\
    1  & -2 & 1      &        &   \\
       & \ddots & \ddots & \ddots & \\
       &        & 1 & -2 & 1 \\
       &        &   & 1 & -2
    \end{bmatrix}.
\]
The semi-discrete equations then take the form of a linear Hamiltonian system
\begin{equation}\label{eq:linear_hamiltonian_system}
    \dot{\mathbf{y}}(t)
    =\begin{bmatrix} \mathbf{0} & I \\[2pt] L_N & \mathbf{0} \end{bmatrix}\mathbf{y}(t).
\end{equation}
A convenient expression for the discrete Hamiltonian is
\[
    \mathrm{H}(\mathbf{y})
    = \frac{1}{2}\mathbf{p}^\top\mathbf{p} + \frac{1}{2 {\delta x}^{2}}\sum_{j=0}^{N-1}(u_{j+1}-u_j)^2
    \;=\; \frac{1}{2}\mathbf{p}^\top\mathbf{p} - \frac{1}{2}\mathbf{u}^\top L_N\mathbf{u},
\]
so that, with the canonical symplectic matrix
\[
\mathcal{J}=\begin{bmatrix} \mathbf{0} & I \\[2pt] -I & \mathbf{0} \end{bmatrix},
\]
one has the Hamiltonian form
\[
\dot{\mathbf{y}}(t) = \mathcal{J}\nabla \mathrm{H}(\mathbf{y}(t)),
\]
which is equivalent to \eqref{eq:linear_hamiltonian_system}.  This symplectic structure motivates the use of structure-preserving time integrators.  Since \eqref{eq:linear_hamiltonian_system} is linear, many efficient solvers are available (e.g., methods based on the matrix exponential.) Here, the linear case serves as a baseline for comparisons with the nonlinear problem treated next.
\begin{figure}[htbp]
    \centering
    \includegraphics[width=\columnwidth]{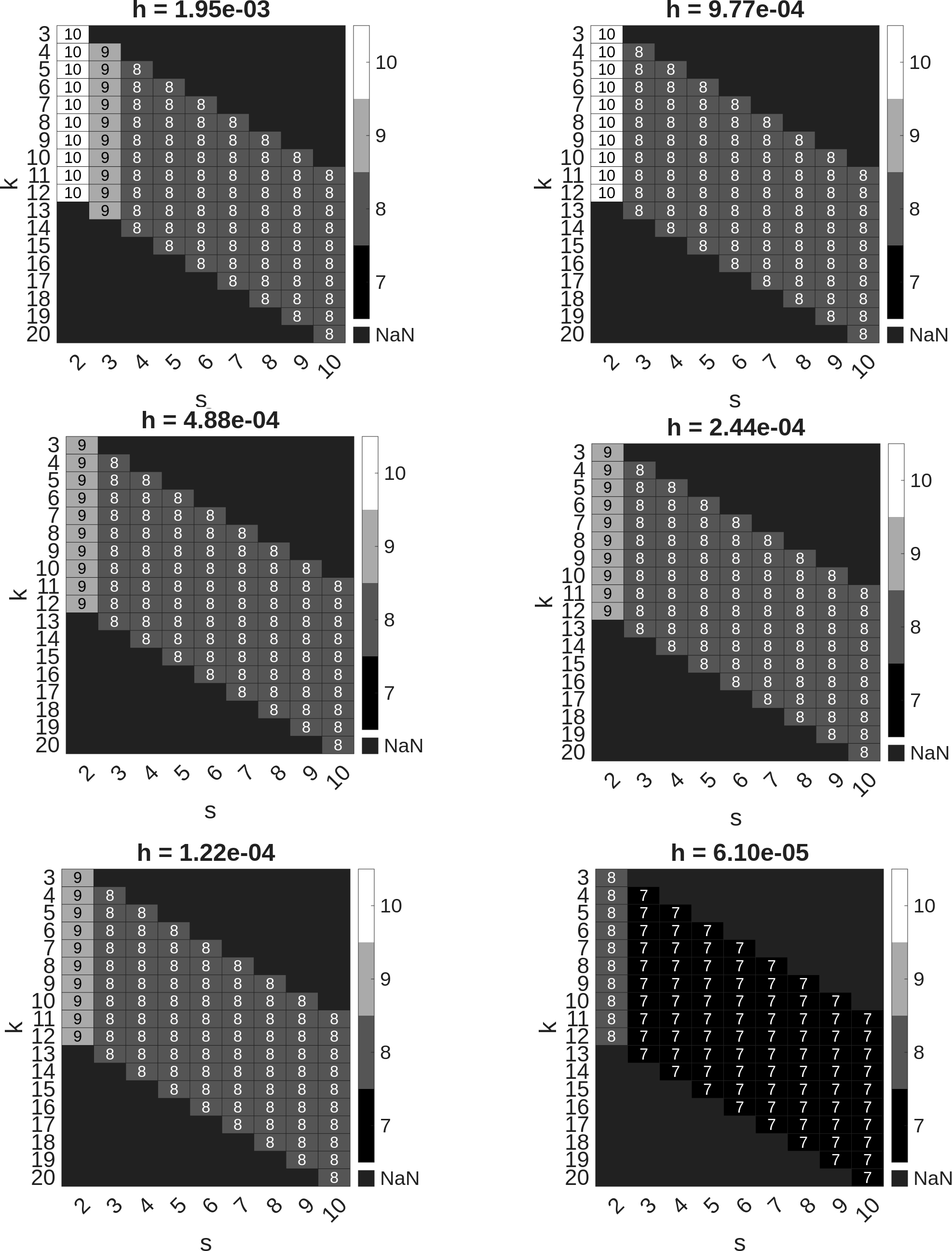}
    \caption{Linear case. We report here the number of iteration of the Krylov projection method for the solution of the Sylvester matrix equation while varying $N = 2^j$, $j=9,\ldots,14$, $h = \nicefrac{1}{N}$, $s = 2,\ldots,10$, and $k = s+1,\ldots,s+10$.}
    \label{fig:linear-case}
\end{figure}
We report in Fig.~\ref{fig:linear-case} the average number of iterations---rounded to the nearest integer---for the solution of the matrix equation~\eqref{eq:sylvester} which is obtained by employing the version of the projection Algorithm~\ref{alg:oneside_sylvester} based on the Extended Krylov method~\cite{MR2318706}, i.e., the overall size of the basis of the Krylov subspace is of the order of twice the number of iterations. For all cases, the inner tolerance for the matrix equation solver is set to $10^{-10}$; the inversion of the $\mathcal{G}$ matrix in this formulation is obtained by computing and storing a sparse LU factorization via UMFPACK~\cite{MR2075981}. We discretize the linear wave equation with $N = 2^j$, $j=9,\ldots,14$, and select for the HBVM($k$,$s$) method a time step $h = \nicefrac{1}{N}$. For the values of $s$ and $k$ we let $s = 2,\ldots,10$ and choose $k = s+1,\ldots,s+10$ hence avoiding the collocation RK case; which we recall corresponds to the choice $s = k$. What we observe, is that the Krylov method behaves robustly with the level of refinement, i.e, the number of iterations does not increase with the decrease of $h$.
\begin{figure}[htbp]
    \centering
    \includegraphics[width=\columnwidth]{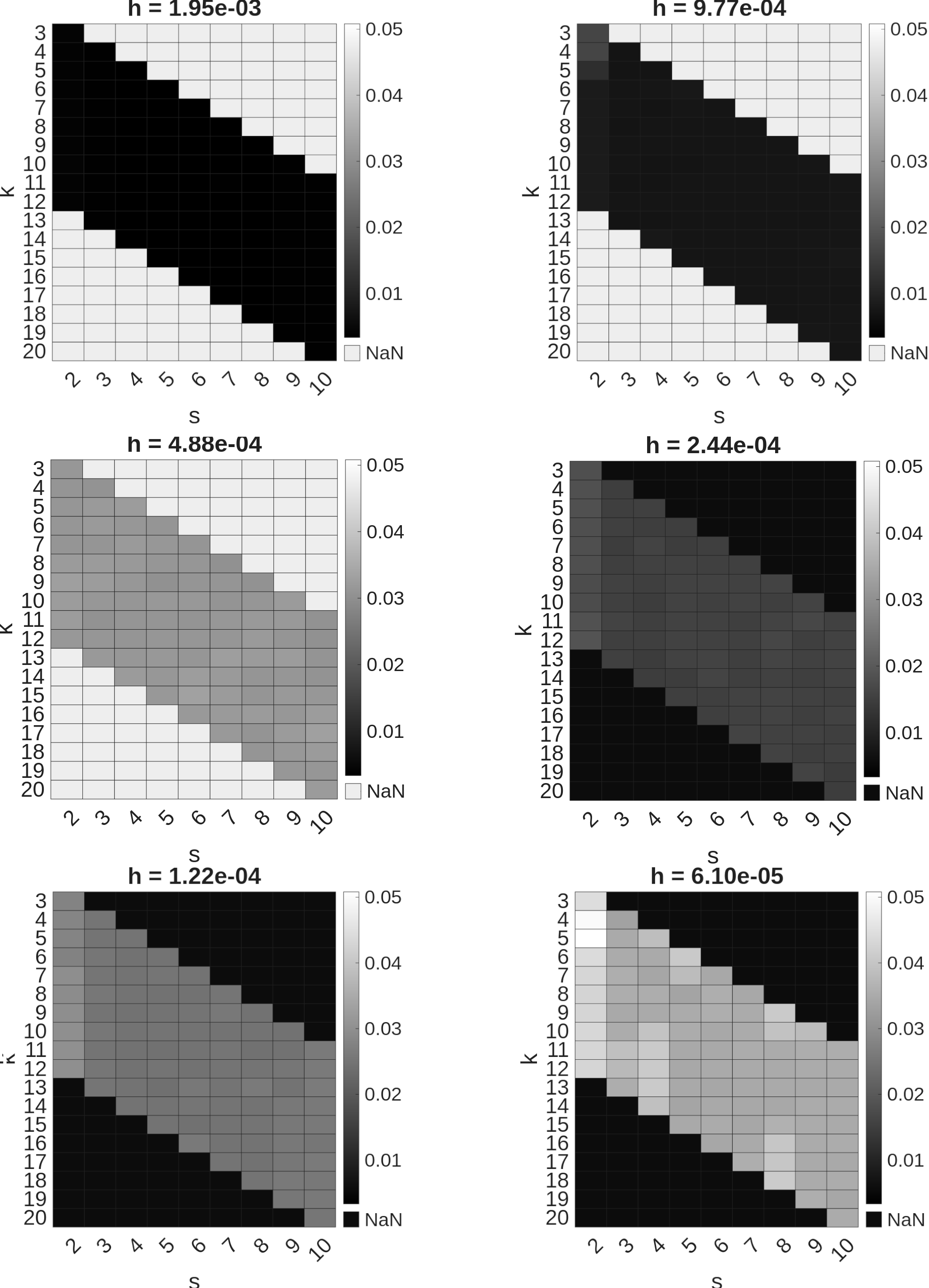}
    \caption{Linear case. We report here the time per step of the overall HBVM($k$,$s$) varying $N = 2^j$, $j=9,\ldots,14$, $h = \nicefrac{1}{N}$,  $s = 2,\ldots,10$, and $k = s+1,\ldots,s+10$.} \label{fig:linear-case-time}
\end{figure}
In Fig.~\ref{fig:linear-case-time} we evaluate the time-per-step of the HBVM($k$,$s$) method while varying the parameters in the same way. This represents a measure of the increase in computational cost given by the mesh-refinement, both in the space and in the time directions. What we observe is the expected behavior, if we double the number of degrees of freedom we are around a doubling in computational time. In this framework, this is the expected scaling, since the overall number of iteration remains stable (see again Fig.~\ref{fig:linear-case}) and the cost-per-iteration is proportional to $N$.

\subsection{A semi-discretized semilinear wave equation}\label{sec:nonlinear_wave}

As a nonlinear extension we consider the semilinear wave equation with a local nonlinear potential, see, e.g.,~\cite{MR4167091}:
\begin{equation}\label{eq:nonlinearwave}
\begin{cases}
    u_{tt}(x,t) = u_{xx}(x,t) + f'(u(x,t)), & (x,t)\in(-L,L)\times[0,\infty),\\[4pt]
    u(x,0)=\psi_0(x),\qquad u_t(x,0)=\psi_1(x), &
\end{cases}
\end{equation}
supplemented by homogeneous Dirichlet boundary conditions
\[
u(-L,t)=0,\qquad u(L,t)=0,\qquad t\ge0.
\]
Introducing the same spatial discretization as above and the state vector \(\mathbf{y}\) as in~\eqref{eq:state-vector}, the semi-discrete system reads
\begin{equation}\label{eq:semidiscrete_nonlinear}
    \dot{\mathbf{y}}(t)
    = \begin{bmatrix}\mathbf{0} & I \\[2pt] L_N & \mathbf{0}\end{bmatrix}\mathbf{y}(t)
      + \begin{pmatrix}\mathbf{0}\\[2pt] -\,f'(\mathbf{u}(t))\end{pmatrix},
\end{equation}
where $f'(\mathbf{u})$ denotes the vector with components $f'(u_j)$.  The corresponding discrete Hamiltonian may be written as
\[
\mathrm{H}(\mathbf{u},\mathbf{p})=\frac{1}{2}\mathbf{p}^\top\mathbf{p}
    + \frac{1}{2 {\delta x}^{2}}\sum_{j=0}^{N-1}(u_{j+1}-u_j)^2 + \sum_{j=1}^{N-1} f(u_j),
\]
so that the semi-discrete dynamics remain Hamiltonian (with the same canonical symplectic form).  

For analysis and for implicit or linearly implicit time integrators in~\eqref{eq:matrix_form_of_simplified_Newton} we need the Jacobian of the right-hand side in~\eqref{eq:semidiscrete_nonlinear}.  Denoting by $\operatorname{diag}(f''(\mathbf{u}))$ the diagonal matrix with entries $f''(u_j)$, the Jacobian takes the block form
\[
    J_F(t,\mathbf{y})
    = \begin{bmatrix}
        \mathbf{0} & I \\[2pt]
        L_N + \operatorname{diag}(f''(\mathbf{u}(t))) & \mathbf{0}
    \end{bmatrix}.
\]
\begin{figure}[htb]
    \centering
    \includegraphics[width=\columnwidth]{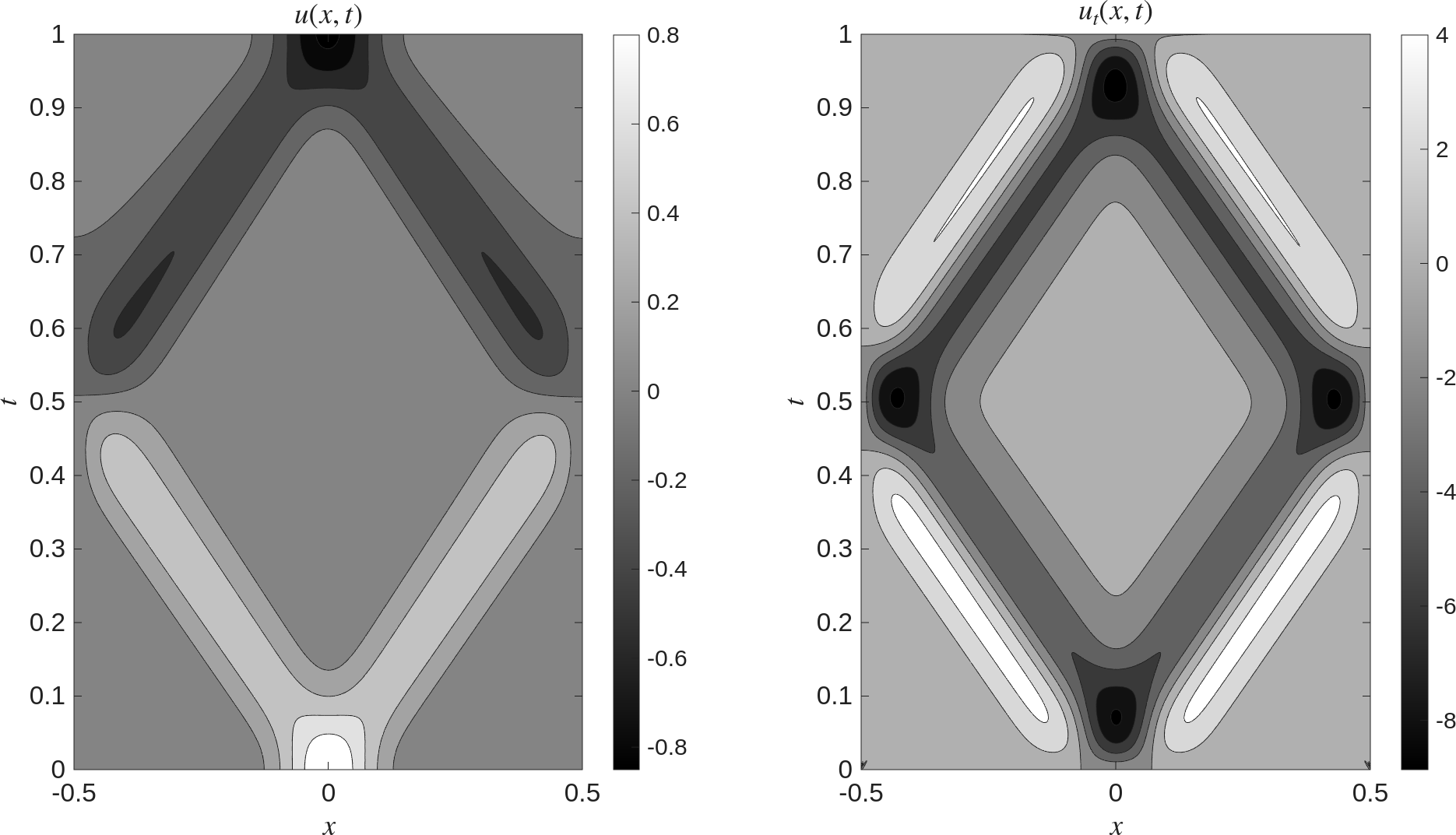}
    
    \caption{Depiction of the numerical solution of~\eqref{eq:nonlinearwave} with $f'(x) = 10 x^2$, $\psi_0(x) = \exp(-100 x^2)$, $\psi_1(x) = 0$.}
    \label{fig:nonlinear_solution}
\end{figure}

To warm-start the Newton iteration for the computation of $\Phi$ at a given time-step, we need a suitable initial choice. Given the form of~\eqref{eq:f=0_in_matrix_form} a suitable choice is to employ the fixed point iteration
\begin{equation}\label{eq:warmup}
\begin{split}
    \Phi^{(0)} = &\; 0, \\
\Phi^{(q+1)}  = &\; \left[ \left. f\left( \mathbf{y}^{(n)} + \overline{h} \Phi^{(q)} \mathcal{I}_s^\top \mathbf{e}_1  \right)  \right| \cdots \left|  f\left( \mathbf{y}^{(n)} + \overline{h} \Phi^{(q)} \mathcal{I}_s^\top \mathbf{e}_{k}  \right) \right. \right] \mathcal{B} \mathcal{W}_s,
\end{split}
\end{equation}
where we start with $\overline{h} = h$, and in case $\|\Phi^{(q+1)} - \Phi^{(q)}\|_F > \|\Phi^{(q)} - \Phi^{(q-1)}\|_F$ we step back and halve $\overline{h} = \nicefrac{\overline{h}}{2}$---indeed, to have guarantees of convergence of the fixed-point iteration this time-step reduction maybe needed. We stop the fixed point iteration when $\|\Phi^{(q+1)} - \Phi^{(q)}\|_F < 10^{-2}$, and, in any case, we never perform more than $20$ fixed point iteration per time-step.

\begin{figure}[htbp]
    \centering
    \input{nonlinearcase1_newtit}

    \includegraphics[width=0.55\columnwidth]{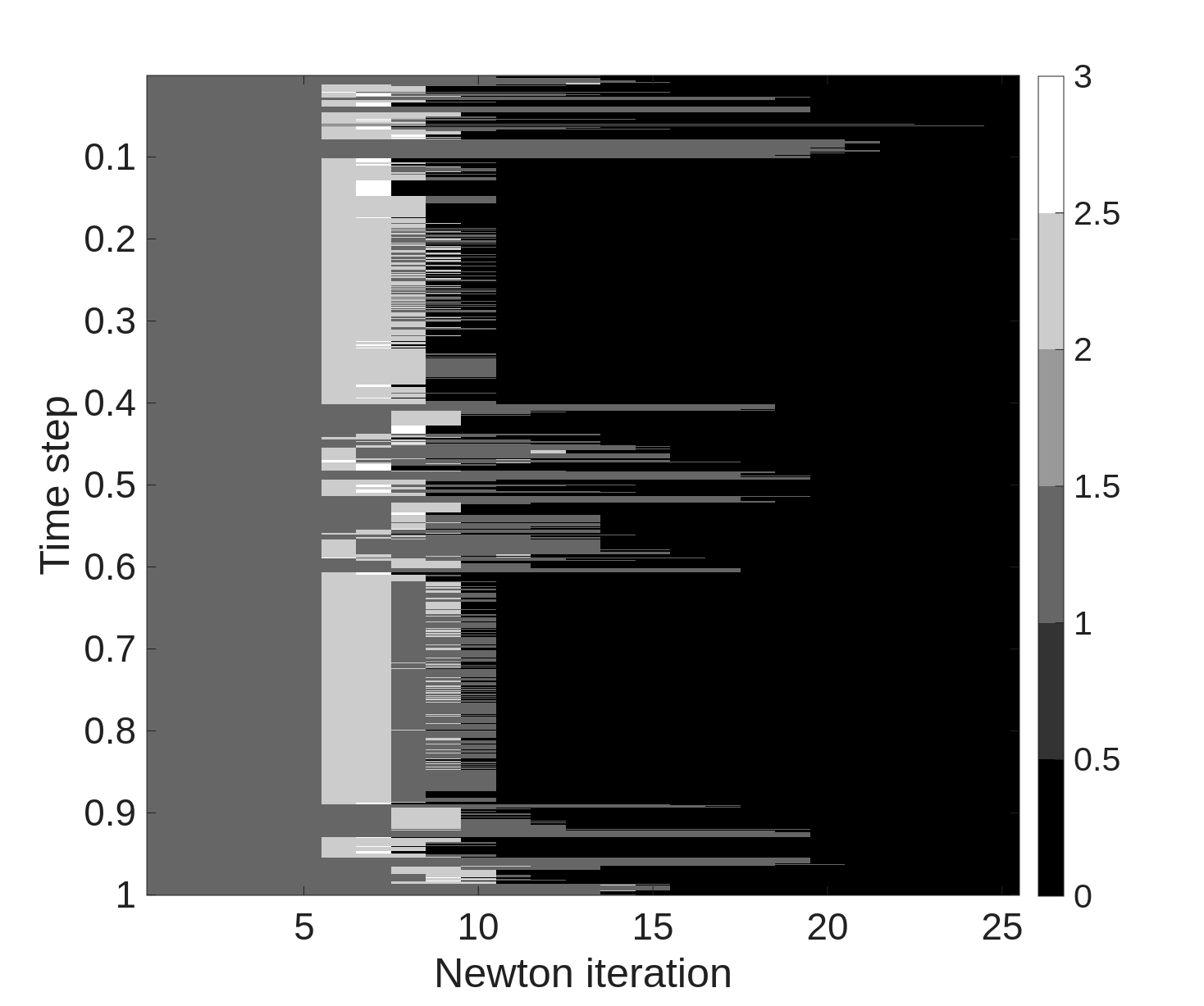}
    \includegraphics[width=0.55\columnwidth]{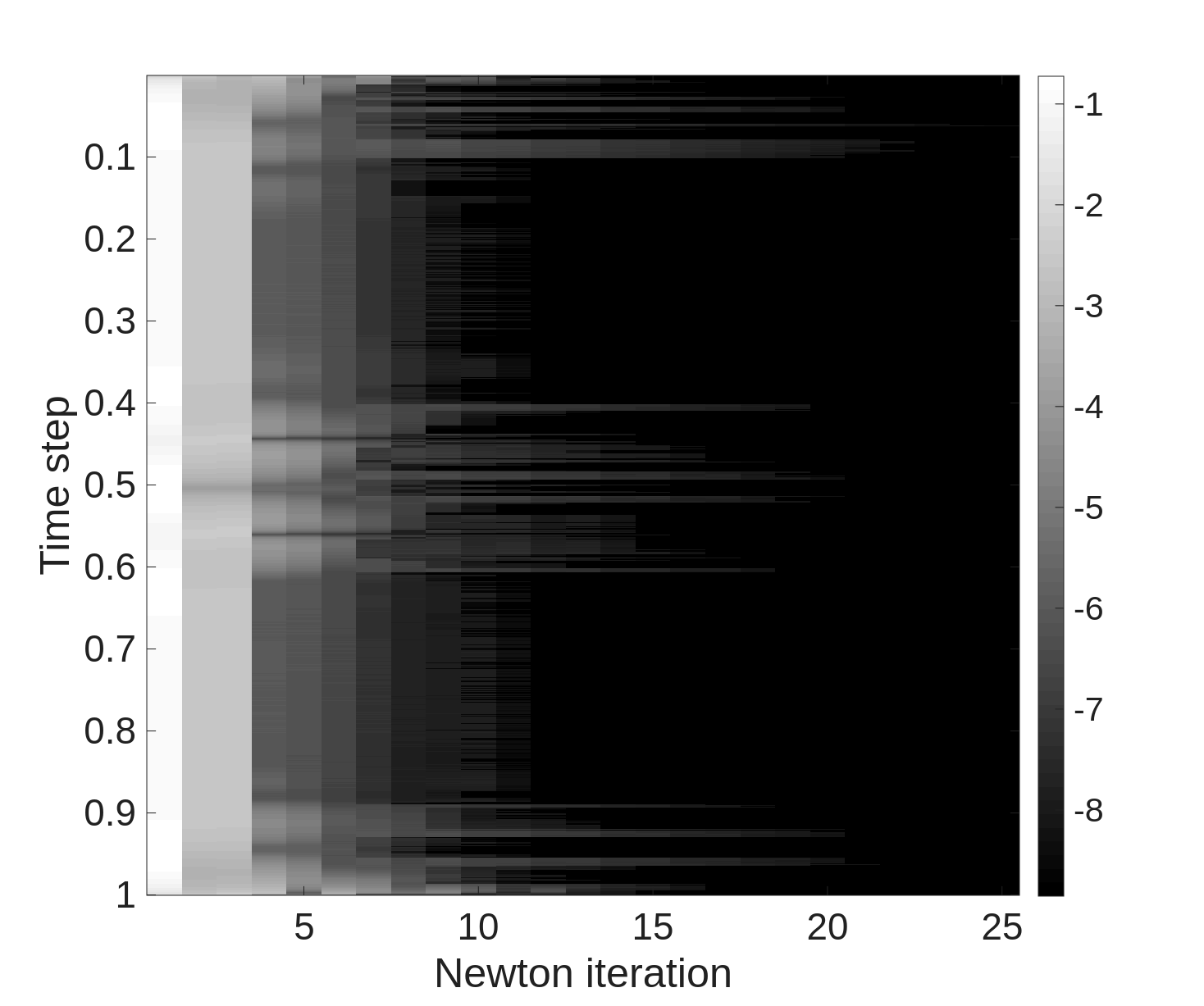}
    
    \caption{Nonlinear test case 1 ($s=2, k=3$). The top figure reports an histogram of the number of Newton iterations per time-step on the left $y$-axis, while on the right there is the residual measured as the norm of $\|\mathcal{F}(\Phi)\|_F$. On the left plot of the second row we have an heatmap with the number of FGMRES iterations preconditioned with the matrix equation solver preconditioner per time step VS Newton iteration; on the right there is the corresponding FGMRES residual, showing the adaptive selection of the tolerance along the Newton iteration.}
    \label{fig:nonlinear1}
\end{figure}

We consider a test case with $f'(u) = 10 u^2$, $N = 1024$, $s = 2$, $k = 3$ and $h = \nicefrac{1}{N}$ for which we report a depiction of the solution in Fig.~\ref{fig:nonlinear_solution}. We report in Fig.~\ref{fig:nonlinear1} the algorithmic performance of the proposed approach.
Specifically, we consider using the complete Newton method with an absolute tolerance of $10^{-8}$ and relative of $10^{-10}$; indeed the actual residual, reported on the right $y$-axis of the topmost panel of Fig.~\ref{fig:nonlinear1}, is between $10^{-9}$ and $10^{-8}$. We set the maximum number of Newton iterations to $100$ per time-step; importantly, the time-step adaptivity mechanism described in Section~\ref{sec:preconditioner} is never triggered throughout the entire integration, indicating robust convergence without the need for step-size reductions. We observe that the number of Newton iterations remains under-control, i.e., it is in between $8$ and $25$ with the majority of steps being nearer to the lower bound. For all time-steps, we never do more than 5 fixed point iterations~\eqref{eq:warmup} to produce an initial guess; in almost all cases we have to reduce one time the $\overline{h}$ time-step for the fixed-point iteration. To solve the linear system with the Jacobian we employ the FGMRES method with a tolerance on the convergence decided automatically by the Newton method---as described in Section~\ref{sec:deciding_tolerances}---and reported in the bottom-right panel of Fig.~\ref{fig:nonlinear1}; as expected the requirement on the tolerance is proportional to the magnitude of the residual in the given Newton iteration. The bottom-left panel of Fig.~\ref{fig:nonlinear1} shows the number of FGMRES iteration needed to achieve such residual, and as we read from the number of iterations the preconditioner~\eqref{eq:matrix_equation_for_Newton} is highly effective. Indeed, we make no more than 3 outer iterations to reach the required accuracy. Note that this is consistent with the good convergence properties of the inner Krylov method for the matrix equation shown in the analysis for the linear case; see again Fig.~\ref{fig:linear-case}.

\begin{figure}[htbp]
    \centering
    \input{nonlinearcase2_newtit}

    \includegraphics[width=0.55\columnwidth]{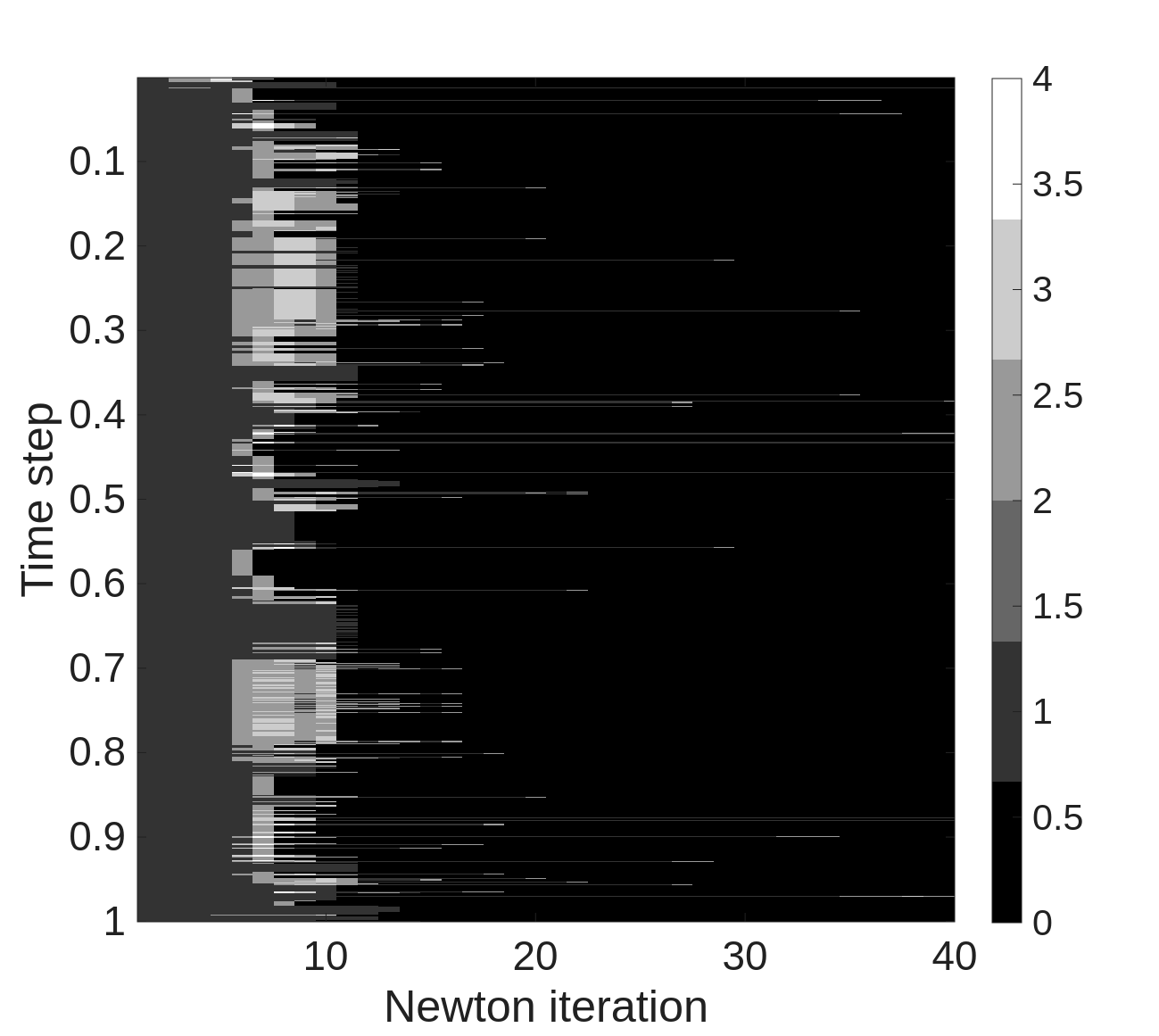}
    \includegraphics[width=0.55\columnwidth]{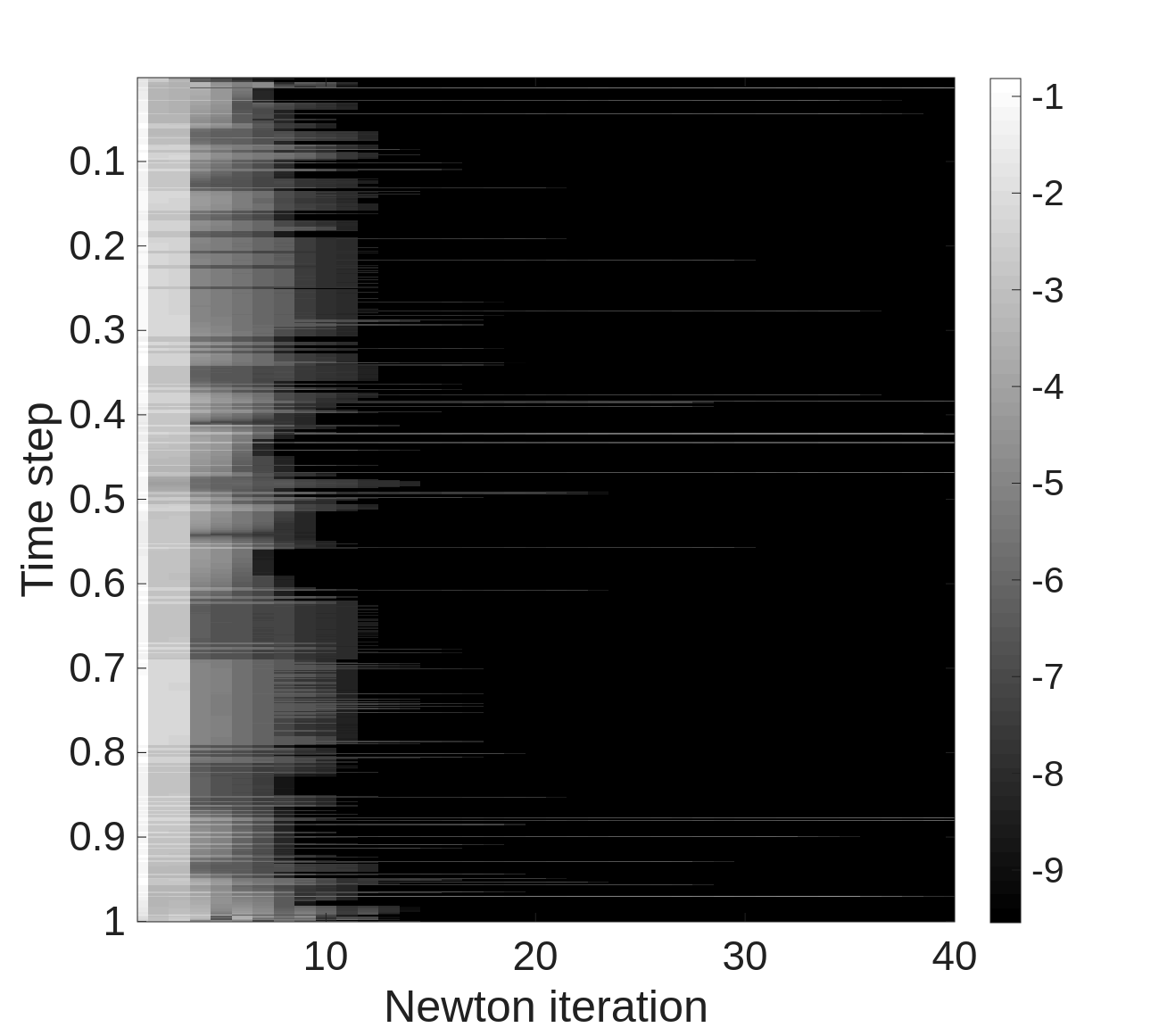}
    
    \caption{Nonlinear test case 2 ($s=3, k=6$). The top figure reports an histogram of the number of Newton iterations per time-step (averaged every 5 time-steps) on the left $y$-axis, while on the right there is the residual measured as the norm of $\|\mathcal{F}(\Phi)\|_F$ for each time-step. On the left plot of the second row we have an heatmap with the number of FGMRES iterations preconditioned with the matrix equation solver preconditioner per time-step VS Newton iteration; on the right there is the corresponding FGMRES residual, showing the adaptive selection of the tolerance along the Newton iteration. In the bottom panels, we have cropped the number of Newton iterations to 40, consistently with the top panel.}
    \label{fig:nonlinear2}
\end{figure}

To look at a different configuration of the HBVM method, we consider also the same test case with $f'(u) = 10 u^2$, $N = 1024$, and $h = \nicefrac{1}{N}$ but selecting $s = 3$, $k = 6$. We report the corresponding results in Fig.~\ref{fig:nonlinear2}. 
For this higher-order configuration, we employ the same absolute and relative tolerances of $10^{-8}$ and $10^{-10}$ respectively for the complete Newton method, resulting in actual residuals in the range of $10^{-9}$ to $10^{-8}$ as shown on the right $y$-axis of the topmost panel of Fig.~\ref{fig:nonlinear2}. We also set the maximum number of Newton iterations to $100$ per time-step; in contrast to the first test case, the time-step adaptivity mechanism is triggered for a few steps during the integration, where the Newton iteration fails to converge within the prescribed limit and the time-step is consequently halved and the computation restarted. Despite the increased values of $s$ and $k$ compared to the first test case, the Newton iteration count remains well-controlled, ranging from approximately $7$ to $40$ iterations per time-step, with most steps requiring between $10$ and $15$ iterations. The warm-start strategy based on the fixed point iteration~\eqref{eq:warmup} continues to be effective, requiring no more than 5 iterations to generate a suitable initial guess for the Newton method. Similarly to the previous case, we occasionally need to halve the time-step $\overline{h}$ during the fixed-point iteration to ensure convergence. The performance of the FGMRES method with the matrix equation solver preconditioner, illustrated in the bottom panels of Fig.~\ref{fig:nonlinear2}, demonstrates that the preconditioner~\eqref{eq:matrix_equation_for_Newton} maintains its effectiveness for this higher-order configuration. The number of outer FGMRES iterations remains modest, typically not exceeding 4 iterations to achieve the required accuracy dictated by the adaptive tolerance strategy shown in the bottom-right panel. This confirms that the proposed preconditioner scales well with the method parameters, maintaining robust convergence properties even when the number of stages $s$ and $k$ increases.

\section{Conclusions and future extensions}
\label{sec:conclusions}

In this work, we have developed efficient solution strategies for Hamiltonian Boundary Value Methods (HBVMs) applied to Hamiltonian systems by exploiting the inherent low-rank structure of their stage equations. For linear Hamiltonian systems, we demonstrated that the stage equations can be reformulated as matrix equations with low-rank right-hand sides, which can be solved efficiently using Krylov projection methods. The numerical experiments on semi-discretized linear wave equations showed that the Krylov solver exhibits robust convergence behavior with iteration count that remains stable across mesh refinement levels, confirming excellent scalability properties.

For nonlinear Hamiltonian systems, we extended the low-rank approach in two complementary ways. First, we showed how the simplified Newton method naturally inherits the low-rank matrix equation structure, enabling efficient iterative solution of the stage equations. Second, and more significantly, we developed a preconditioner for the complete Newton--Krylov framework based on solving an approximate matrix equation at each linear iteration. This preconditioner, implemented within a FGMRES solver with adaptive forcing terms, proved highly effective in our numerical experiments. The tests on semi-discretized semilinear wave equations demonstrated that only a few outer FGMRES iterations are required to achieve the prescribed accuracy, even for higher-order HBVM configurations. We further adopted an adaptive time-stepping strategy that dynamically adjusts the step size to ensure Newton convergence, providing additional robustness for challenging nonlinear problems.

The construction discussed here should be applicable also to initial value problems involving Caputo fractional derivatives, for which HBVM methods have been extended in~\cite{MR4892878}; therefore, we plan to extend the low-rank procedure to this case as well. Furthermore, from an implementation standpoint, the proposed approach becomes particularly relevant when the matrices involved are so large that storing them on a single machine is impractical. In such cases, developing a distributed implementation of the matrix equation solver is of great interest. Looking ahead, we plan to extend the current implementation using the \texttt{PSCToolkit}~\cite{DAmbra2023} library, which would enable both distributed computing and GPU acceleration for efficiently solving large-scale problems. Another promising direction is the application of these techniques to other structure-preserving integrators with similar low-rank properties, potentially extending the benefits of this approach to a broader class of geometric numerical methods. Finally, in future work we also plan to investigate specialized pole–selection strategies tailored to the spectral properties of the matrices arising in our formulation for more complex Hamiltonian PDEs, with the aim of further accelerating convergence in rational Krylov–based solvers.

\begin{acknowledgement}
FD acknowledges the MUR Excellence Department Project awarded to the Department of Mathematics, University of Pisa, CUP I57G22000700001.  MM acknowledges the MUR Excellence Department Project MatMod@TOV awarded to the Department of Mathematics, University of Rome Tor Vergata, CUP E83C23000330006. The research of FD was  partially granted by the Italian Ministry of University and Research (MUR) through the PRIN 2022 ``MOLE: Manifold constrained Optimization and LEarning'',  code: 2022ZK5ME7 MUR D.D. financing decree n. 20428 of November 6th, 2024 (CUP B53C24006410006). The research of MM was  partially granted by the Italian Ministry of University and Research (MUR) through the PRIN 2022-PNRR ``MATHPROCULT:
MATHematical tools for predictive maintenance and PROtection of CULTtural heritage'' (CUP J53D23015940001). Both authors are members of INdAM-GNCS and have been partially financed by the INdAM-GNCS Project CUP E53C24001950001.

We thank the anonymous Referees for the careful reading of the manuscript, and their suggestions which have improved it.
\end{acknowledgement}
\ethics{Competing Interests}{
The authors have no conflicts of interest to declare that are relevant to the content of this chapter.}

\bibliographystyle{spmpsci}
\bibliography{bibliography}

\end{document}